\documentclass[a4paper,11pt]{article}

\usepackage[T1]{fontenc}
\usepackage{lmodern}

\usepackage{dsfont}

\usepackage[latin1]{inputenc}
\usepackage{amsmath}
\usepackage{amsthm}
\usepackage{amssymb}
\usepackage{mathrsfs}
\usepackage{graphicx}
\usepackage[all]{xy}
\usepackage{hyperref}

\usepackage{stmaryrd}
\usepackage{caption}

\usepackage{abstract} 

\newtheorem{thm}{Theorem}[section]
\newtheorem{cor}[thm]{Corollary}
\newtheorem{claim}[thm]{Claim}

\newtheorem{lemma}[thm]{Lemma}
\newtheorem{prop}[thm]{Proposition}

\theoremstyle{definition}
\newtheorem{definition}[thm]{Definition}

\title{A note on torsion subgroups of groups acting on finite-dimensional CAT(0) cube complexes}
\date{\today}
\author{Anthony Genevois}

\begin{document}

\maketitle

\begin{abstract}
In this article, we state and prove a general criterion which prevent some groups from acting properly on finite-dimensional CAT(0) cube complexes. As an application, we show that, for every non-trivial finite group $F$, the lamplighter group $F \wr \mathbb{F}_2$ over a free group does not act properly on a finite-dimensional CAT(0) cube complex (although it acts properly on a infinite-dimensional CAT(0) cube complex). We also deduce from this general criterion that, roughly speaking, given a group $G$ acting on a CAT(0) cube complex of finite dimension and an infinite torsion subgroup $L \leq G$, either the normaliser $N_G(L)$ is close to be free abelian or, for every $k \geq 1$, $N_G(L)$ contains a non-abelian free subgroup commuting with a subgroup of $L$ of size $\geq k$. 
\end{abstract}

\tableofcontents

\section{Introduction}

In the last decades, it has been proved that finding a proper action of a given group on a CAT(0) cube complex provides interesting information on the group. For instance, if such an action exists, then the group has to satisfy the Haagerup property \cite{Haagerup}, its polycyclic subgroups must be virtually abelian \cite{CornulierCommensurated, CFTT}, and as soon as our group is finitely generated its free abelian subgroups are necessarily undistorted \cite{CornulierCommensurated, WoodhouseFlatTorus}. But if it turns out that the cube complex can be chosen finite-dimensional, then some additional information is obtained, including weak amenability \cite{Mizura, WeakAmenability}, RD property \cite{RDccc} and Tits alternative \cite{TitsCCC, CapraceSageev}. 

\medskip \noindent
Consequently, given a group acting properly on a CAT(0) cube complex, a natural question is: can the cube complex be finite-dimensional? An easy obstruction is given the following observation: if a group acts properly on a CAT(0) cube complex of dimension $d$, then a free abelian subgroup must have rank at most $d$. Therefore, a group containing a free abelian group of arbitrary large rank cannot act properly on a finite-dimensional CAT(0) cube complex. Such groups include for instance Thompson-like groups. But otherwise, if there exists a bound on the rank of free abelian subgroups, the question is delicate.

\medskip \noindent
The main result of this article provides a general criterion preventing some groups from acting properly on finite-dimensional CAT(0) cube complexes. More precisely, we prove:

\begin{thm}\label{intro:MainCriterion}
Let $\mathcal{C}$ be a collection of finitely generated groups. Assume that:
\begin{itemize}
	\item for every group $G \in \mathcal{C}$ and every finite-index subgroup $H \leq G$, there exists a subgroup of $H$ which belongs to $\mathcal{C}$;
	\item for every group $G \in \mathcal{C}$, there exists a subgroup of the commutator subgroup $[G,G]$ which belongs to $\mathcal{C}$;
	\item every group of $\mathcal{C}$ contains a non-abelian free subgroup and a $\delta$-normal subgroup which is (locally finite)-by-(free abelian).
\end{itemize}
Then no group of $\mathcal{C}$ acts properly on a finite-dimensional CAT(0) cube complex. 
\end{thm}

\noindent
Given a group $G$ and a subgroup $H \leq G$, one says that $H$ is \emph{$\delta$-normal} if, for every finite collection $g_1, \ldots, g_n \in G$, the intersection $\bigcap\limits_{i=1}^n g_iHg_i^{-1}$ is infinite.

\medskip \noindent
In the specific case where the collection $\mathcal{C}$ is reduced to a single finitely generated group, one obtains the following statement:

\begin{cor}
Let $G$ be a finitely generated group. Assume that:
\begin{itemize}
	\item the commutator subgroup $[G,G]$ and every finite-index subgroup of $G$ contains a copy of $G$;
	\item $G$ contains a non-abelian free subgroup and a $\delta$-normal subgroup which is (locally finite)-by-(free abelian).
\end{itemize}
Then $G$ does not act properly on a finite-dimensional CAT(0) cube complex.
\end{cor}

\noindent
For instance, this criterion applies to \emph{wreath products}, namely:

\begin{prop}\label{intro:Wreath}
For every non-trivial finite group $F$, the wreath product $F \wr \mathbb{F}_2$ does not act properly on a finite-dimensional CAT(0) cube complex.
\end{prop}

\noindent
Recall that, given two groups $A$ and $B$, the wreath product $A \wr B$ is defined as the semi-direct product 
$$\left( \bigoplus\limits_{g \in B} A \right) \rtimes B,$$
where $B$ acts on the direct sum by permuting the coordinates. It has been proved in \cite{Haagerupwreath} (see also \cite{LampMedian}) that acting properly on a CAT(0) cube complex is a property which is stable under wreath products. However, it is not clear when such a complex can be taken finite-dimensional. For instance, it is proved in \cite[Proposition 9.33]{Qm} that, for every finite group $F$ and every $n \geq 1$, the wreath product $F \wr \mathbb{Z}^n$ acts properly on a CAT(0) cube complex of dimension $2n$. On the other hand, as observed in \cite{CSVparticulier}, the wreath product $F \wr \mathbb{F}_2$ does not act \emph{metrically properly} on a finite-dimensional CAT(0) cube complex. So Proposition \ref{intro:Wreath} improves slightly this observation by replacing ``metrically proper'' with ``proper''. 

\medskip \noindent
The argument of \cite{CSVparticulier} is indirect, based on weak amenability. The main motivation of the article was to understand cubically what prevents $F \wr \mathbb{F}_2$ from acting properly on a finite-dimensional CAT(0) cube complex. 

\medskip \noindent
Another interesting consequence of Theorem \ref{intro:MainCriterion} is more theoretical. It provides information on normalisers of infinite torsion subgroups of groups acting properly on finite-dimensional CAT(0) cube complexes. Roughly speaking, given a group $G$ acting on a CAT(0) cube complex of finite dimension and an infinite torsion subgroup $L \leq G$, one shows that either the normaliser $N_G(L)$ is close to be free abelian or, for every $k \geq 1$, $N_G(L)$ contains a non-abelian free subgroup commuting with a subgroup of $L$ of size $\geq k$. More precisely:

\begin{thm}\label{intro:torsionsub}
Let $G$ be a group acting properly on a finite-dimensional CAT(0) cube complex $X$. Assume that $G$ contains an infinite torsion subgroup $L$. Either the normaliser $N_G(L)$ contains a finite-index subgroup which is (locally finite)-by-$\mathbb{Z}^n$ for some $n \leq \dim(X)$; or, for every $k \geq 1$, $N_G(L)$ contains a contains a non-abelian free subgroup centralising a subgroup of $L$ of cardinality $\geq k$.
\end{thm}

\noindent
Notice that Proposition \ref{intro:Wreath} also follows from this statement. As $F \wr \mathbb{F}_2$ acts properly on an infinite-dimensional CAT(0) cube complex, it follows that the conclusion of Theorem~\ref{intro:torsionsub} does not hold if we drop the hypothesis of finite dimension.

\paragraph{Acknowledgments.} This work was supported by a public grant as part of the Fondation Math\'ematique Jacques Hadamard.

\section{Preliminaries}

\noindent
A \textit{cube complex} is a CW complex constructed by gluing together cubes of arbitrary (finite) dimension by isometries along their faces. It is \textit{nonpositively curved} if the link of any of its vertices is a simplicial \textit{flag} complex (ie., $n+1$ vertices span a $n$-simplex if and only if they are pairwise adjacent), and \textit{CAT(0)} if it is nonpositively curved and simply-connected. See \cite[page 111]{MR1744486} for more information.

\begin{figure}
\begin{center}
\includegraphics[trim={0 13cm 10cm 0},clip,scale=0.45]{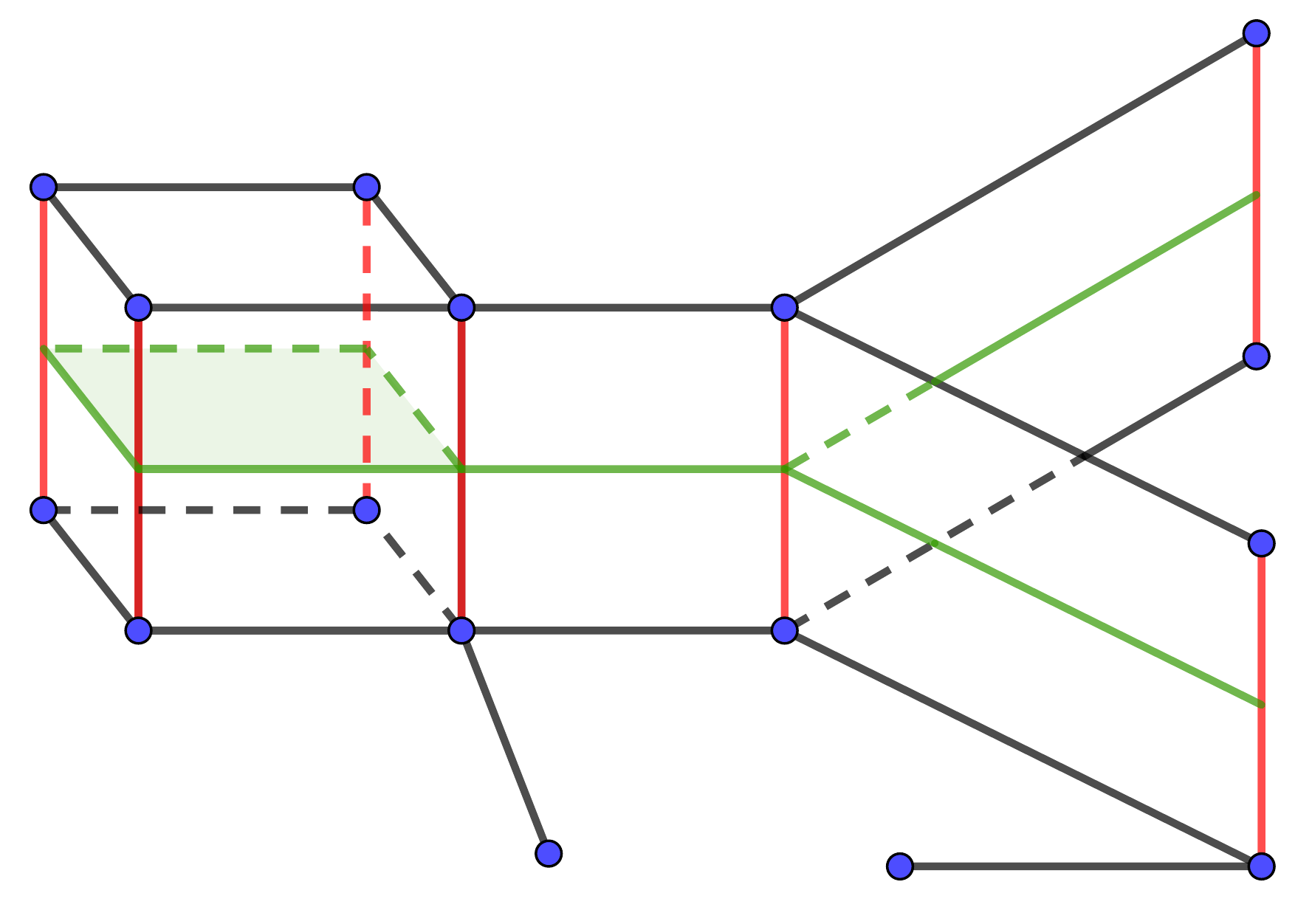}
\caption{A hyperplane (in red) and the associated union of midcubes (in green).}
\label{figure27}
\end{center}
\end{figure}

\medskip \noindent
Fundamental tools when studying CAT(0) cube complexes are \emph{hyperplanes}. Formally, a \textit{hyperplane} $J$ is an equivalence class of edges with respect to the transitive closure of the relation identifying two parallel edges of a square. Geometrically, a hyperplane $J$ is rather thought of as the union of the \textit{midcubes} transverse to the edges belonging to $J$ (sometimes referred to as its \emph{geometric realisation}). See Figure \ref{figure27}. 

\medskip \noindent
There exist several metrics naturally defined on a CAT(0) cube complex. For instance, one can extend in a standard way the Euclidean metrics defined on each cube to a global length metric, and the distance one obtains in this way turns out to be CAT(0). However, in this article, we are mainly interested in the graph metric defined on the one-skeleton of the cube complex, referred to as its \emph{combinatorial metric}. Unless specified otherwise, we will always identify a CAT(0) cube complex with its one-skeleton, thought of as a collection of vertices endowed with a relation of adjacency. 

\medskip \noindent
The following theorem will be often used along the article without mentioning it.

\begin{thm}\emph{\cite{MR1347406}}
Let $X$ be a CAT(0) cube complex.
\begin{itemize}
	\item If $J$ is a hyperplane of $X$, the graph $X \backslash \backslash J$ obtained from $X$ by removing the (interiors of the) edges of $J$ contains two connected components. They are convex subgraphs of $X$, referred to as the \emph{halfspaces} delimited by $J$.
	\item For every vertices $x, y \in X$, the distance between $x$ and $y$ coincides with the number of hyperplanes separating them.
\end{itemize}
\end{thm}

\subsection{Wallspaces, cubical quotients}

\noindent
Given a set $X$, a \emph{wall} $\{A,B\}$ is a partition of $X$ into two non-empty subsets $A,B$, referred to as \emph{halfspaces}. Two points of $X$ are \emph{separated} by a wall if they belong to two distinct subsets of the partition. 

\medskip \noindent
A \emph{wallspace} $(X, \mathcal{W})$ is the data of a set $X$ and a collection of walls $\mathcal{W}$ such that any two points are separated by only finitely many walls. Such a space is naturally endowed with the pseudo-metric
$$d : (x,y) \mapsto \text{number of walls separating $x$ and $y$}.$$
As shown in \cite{ChatterjiNiblo, NicaCubulation}, there is a natural CAT(0) cube complex associated to any wallspace. More precisely, given a wallspace $(X, \mathcal{W})$, define an \emph{orientation} $\sigma$ as a collection of halfspaces such that:
\begin{itemize}
	\item for every $\{A,B\} \in \mathcal{W}$, $\sigma$ contains exactly one subset among $\{A,B\}$;
	\item if $A$ and $B$ are two halfspaces satisfying $A \subset B$, then $A \in \sigma$ implies $B \in \sigma$.
\end{itemize}
Roughly speaking, an orientation is a coherent choice of a halfspace in each wall. As an example, if $x \in X$, then the set of halfspaces containing $x$ defines an orientation. Such an orientation is referred to as a \emph{principal orientation}. Notice that, because any two points of $X$ are separated by only finitely many walls, two principal orientations are always \emph{commensurable}, ie., their symmetric difference is finite.

\begin{figure}
\begin{center}
\includegraphics[trim={0 18cm 27cm 0},clip,scale=0.45]{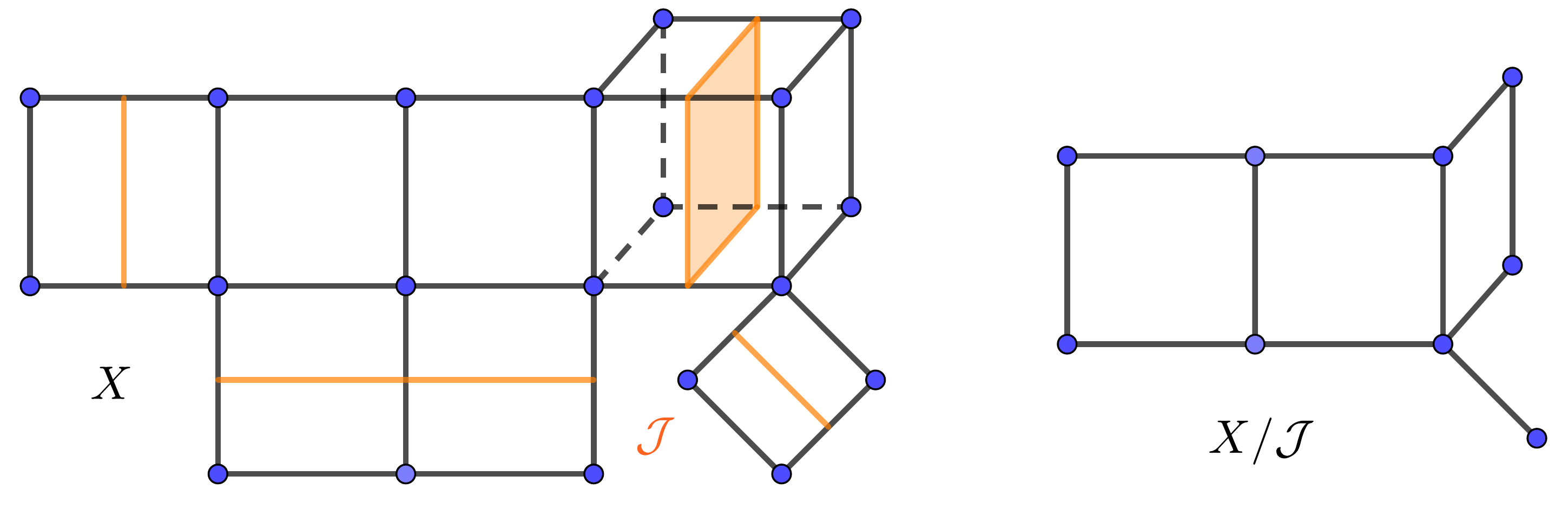}
\caption{A CAT(0) cube complex and one of its cubical quotients.}
\label{quotient}
\end{center}
\end{figure}

\medskip \noindent
The \emph{cubulation} of $(X, \mathcal{W})$ is the cube complex
\begin{itemize}
	\item whose vertices are the orientations within the commensurability class of principal orientations;
	\item whose edges link two orientations if their symmetric difference has cardinality two;
	\item whose $n$-cubes fill in all the subgraphs isomorphic to one-skeleta of $n$-cubes.
\end{itemize}
See Figure \ref{quotient} for an example.

\begin{definition}
Let $X$ be a CAT(0) cube complex and $\mathcal{J}$ be a collection of hyperplanes. Let $\mathcal{W}(\mathcal{J})$ denote the set of partitions of $X$ induced by the hyperplanes of $\mathcal{J}$. The \emph{cubical quotient} $X/ \mathcal{J}$ of $X$ by $\mathcal{J}$ is the cubulation of the wallspace $(X, \mathcal{W}(\mathcal{J}))$. 
\end{definition}

\noindent
See Figure \ref{quotient} for an example. It can be shown that $X/ \mathcal{J}$ can be obtained from $X$ in the following way. Given a hyperplane $J \in \mathcal{J}$, cut $X$ along $J$ to obtain $X \backslash \backslash J$. Each component of $X \backslash \backslash J$ contains a component of $N(J) \backslash \backslash J= N_1 \sqcup N_2$. Notice that $N_1$ and $N_2$ are naturally isometric: associate to each vertex of $N_1$ the vertex of $N_2$ which is adjacent to it in $N(J)$. Now, glue the two components of $X \backslash \backslash J$ together by identifying $N_1$ and $N_2$. The cube complex obtained is still CAT(0) and its set of hyperplanes corresponds naturally to $\mathcal{H}(X) \backslash \{J\}$ if $\mathcal{H}(X)$ denotes the set of hyperplanes of $X$. Thus, the same construction can repeated with a hyperplane of $\mathcal{J} \backslash \{J\}$, and so on. The cube complex which is finally obtained is the cubical quotient $X/ \mathcal{J}$.

\medskip \noindent
Notice that a quotient map is naturally associated to a cubical quotient, namely:
$$\left\{ \begin{array}{ccc} X & \to & X/ \mathcal{J} \\ x & \mapsto & \text{principal orientation defined by $x$} \end{array} \right..$$
The next lemma relates the distance between two vertices of $X$ to the distance between their images in $X/ \mathcal{J}$.

\begin{lemma}\label{lem:DistQuotient}
Let $X$ be a CAT(0) cube complex and $\mathcal{J}$ a collection of hyperplanes. If $\pi : X \to X/ \mathcal{J}$ denotes the canonical map, then
$$d_{X/ \mathcal{J}} \left( \pi(x), \pi(y) \right) = \# \left( \mathcal{W}(x,y) \backslash \mathcal{J} \right)$$
for every vertices $x,y \in X$, where the set $\mathcal{W}(x,y)$ denotes the collection of the hyperplanes of $X$ separating $x$ and $y$. \qed
\end{lemma}

\subsection{Roller boundary}

\noindent
Let $X$ be a CAT(0) cube complex. An \emph{orientation} of $X$ is an orientation of the wallspace $(X, \mathcal{W}(\mathcal{J}))$, as defined in the previous section, where $\mathcal{J}$ is the set of all the hyperplanes of $X$. The \emph{Roller compactification} $\overline{X}$ of $X$ is the set of the orientations of $X$. Usually, we identify $X$ with the image of the embedding
$$\left\{ \begin{array}{ccc} X & \to & \overline{X} \\ x & \mapsto & \text{principal orientation defined by $x$} \end{array} \right.$$
and we define the \emph{Roller boundary} of $X$ by $\mathfrak{R}X:= \overline{X} \backslash X$. 

\begin{figure}
\begin{center}
\includegraphics[trim={0 12cm 41cm 0},clip,scale=0.45]{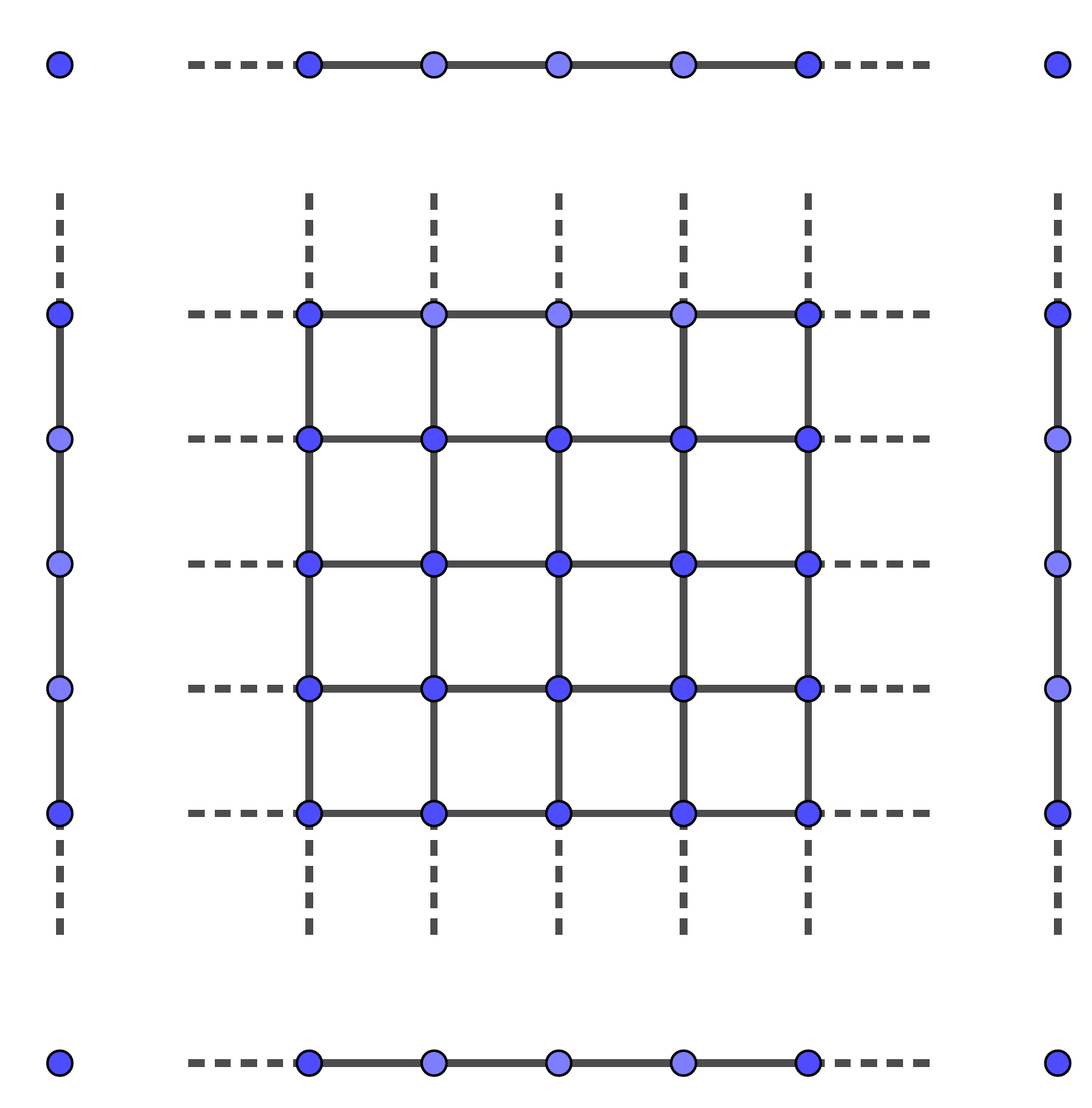}
\caption{Roller compactification of $\mathbb{R}^2$. It contains nine cubical components.}
\label{Roller}
\end{center}
\end{figure}

\medskip \noindent
The Roller compactification is also naturally a cube complex. Indeed, if we declare that two orientations are linked by an edge if their symmetric difference has cardinality two and if we declare that any subgraph isomorphic to the one-skeleton of an $n$-cube is filled in by an $n$-cube for every $n \geq 2$, then $\overline{X}$ is a disjoint union of CAT(0) cube complexes. Each such component is referred to as a \emph{cubical component} of $\overline{X}$. See Figure~\ref{Roller} for an example. Notice that the distance (possibly infinite) between two vertices of $\overline{X}$ coincides with the number of hyperplanes which separate them, if we say that a hyperplane $J$ separates two orientations when they associate distinct halfspaces to $J$. In the sequel, we also say that an orientation $\sigma$ belongs to a halfspace $D$ if $\sigma(J)=D$ where $J$ is the hyperplane bounding $D$.

\medskip \noindent
An elementary but useful observation is that cubical components in the Roller boundary have dimensions smaller than the cube complex. More precisely:

\begin{lemma}
Let $X$ be a finite-dimensional CAT(0) cube complex. For every cubical component $Y$ of $\mathfrak{R}X$, the inequality $\dim(Y)<\dim(X)$ holds.
\end{lemma}

\noindent
Let us conclude this subsection with two remarks. First, the Roller compactification $\overline{X}$ of a CAT(0) cube complex can be endowed with a ternary operation $\mu : \overline{X} \times \overline{X} \times \overline{X} \to \overline{X}$ which extends the median operation in $X$. Namely, 

\begin{definition}
Let $X$ be a CAT(0) cube complex and $\sigma_1,\sigma_2,\sigma_3$ three orientations. The \emph{median point} $\mu(\sigma_1,\sigma_2,\sigma_3)$ is the orientation corresponding to the set of halfspaces which belong to at least two orientations among $\{\sigma_1,\sigma_2,\sigma_3\}$.
\end{definition}

\noindent
Next, if $X$ is a CAT(0) cube complex and $\mathcal{J}$ a collection of hyperplanes, it is worth noticing that the quotient map $\pi : X \to X/ \mathcal{J}$ extends to the Roller compactifications as
$$\overline{\pi} : \left\{ \begin{array}{ccc} \overline{X} & \to & \overline{X/ \mathcal{J}} \\ \sigma & \mapsto & \sigma_{| \mathcal{W} \backslash \mathcal{J}} \end{array} \right.,$$
where $\mathcal{W}$ denotes the collection of all the hyperplanes of $X$.

\subsection{Busemann morphisms}

\noindent
An interesting observation is that, if $T$ is a simplicial tree and $\alpha \in \partial T$ a point at infinity, then the Busemann function $\beta_\alpha$ turns out to define a morphism $\mathrm{stab}(\alpha) \to \mathbb{Z}$. This morphism can be described in the following way. Fixing a basepoint $x \in T$ and an isometry $g \in \mathrm{stab}(\alpha)$, notice that $g \cdot [x,\alpha) \cap [x, \alpha)$ contains a infinite subray. So $g$ acts like a translation on a ray pointing to $\alpha$. The value of $\beta_\alpha(g)$ is precisely the length of this translation (positive it is directed to $\alpha$, and negative otherwise). Moreover, it follows from this description that the kernel of the Busemann morphism coincides with the set of elliptic isometries of $\mathrm{stab}(\alpha)$.

\medskip \noindent
Interestingly, this construction can be extended to finite-dimensional CAT(0) cube complexes, namely:

\begin{thm}\label{thm:Busemann}
Let $X$ be a finite-dimensional CAT(0) cube complex and $\alpha \in \mathfrak{R}X$ a point at infinity. There exist a subgroup $\mathrm{stab}_0(\alpha)$ of $\mathrm{stab}(\alpha) \leq \mathrm{Isom}(X)$ of index at most $\dim(X)!$ and a morphism $\beta : \mathrm{stab}_0(\alpha) \to \mathbb{Z}^n$, where $n \leq \dim(X)$, such that 
$$\mathrm{ker}(\beta)= \{ g \in \mathrm{stab}_0(\alpha) \mid \text{$g$ is $X$-elliptic} \}.$$
Consequently, $\mathrm{stab}(\alpha)$ is virtually (locally $X$-elliptic)-by-(free abelian).
\end{thm}

\noindent
We refer to the appendix of \cite{CFI} and to \cite[Theorem 2.12]{CFTT} for more details.

\subsection{Strongly contracting isometries}

\noindent
In our proof of Theorem \ref{intro:MainCriterion}, the notion of \emph{strongly contracting isometries} will play a fundamental role.

\begin{definition}
Let $X$ be a CAT(0) cube complex. Two hyperplanes $J$ and $H$ are \emph{strongly separated} if no hyperplane is transverse to both $J$ and $H$. An isometry $g \in \mathrm{Isom}(X)$ is \emph{strongly contracting} if there exists two halfspaces $A \subset B$ bounded by two strongly separated hyperplanes such that $g \cdot B \subsetneq A$.
\end{definition}

\noindent
As noticed by \cite[Lemma 6.2]{CapraceSageev}, strongly contracting isometries turn out to define rank-one isometries (or equivalently, contracting isometries). 

\medskip \noindent
In the sequel, the following lemma will be also useful:

\begin{lemma}\label{lem:ContractingPoint}
Let $X$ be a CAT(0) cube complex and $(D_i)$ a decreasing sequence of halfspaces bounded by pairwise strongly separated hyperplanes. Then there exists exactly one point of $\mathfrak{R}X$ which belongs to $D_i$ for every $i$. 
\end{lemma}

\begin{proof}
The orientation $\sigma$ sending each hyperplane $J$ to the halfspace $D$ bounded by $J$ which contains all but finitely many $D_i$ defines a point of $\mathfrak{R}X$ which belongs to $D_i$ for every $i$. 

\medskip \noindent
Next, assume for contradiction that there exist two points $\alpha,\beta \in \mathfrak{R}X$ which belong to $D_i$ for every $i$. Let $J$ be a hyperplane separating $\alpha$ and $\beta$, and fix two vertices $x,y \in X$ separated by $J$. Without loss of generality, suppose that $x$ and $\alpha$ belong to the same halfspace delimited by $J$. If $a \in X$ is a fixed vertex which does not belong to $D_0$, it follows from the fact that there exist only finitely many hyperplanes separating $a$ from either $x$ or $y$ that there exists some $k$ such that $x$ and $y$ do not belong to $D_i$ for every $i \geq k$. Consequently, $J$ separates $\{ \alpha,x \}$ and $\{ \beta, y\}$ and, for every $i \geq k$, $D_i$ contains $\alpha$ and $\beta$ but not $x$ nor $y$. We conclude that $J$ must be transverse to the hyperplane bounding $D_i$ for every $i \geq k$, which is impossible. 
\end{proof}

\noindent
For convenience, we introduce the following definition for future use:

\begin{definition}
Let $X$ be a CAT(0) cube complex. A point $\alpha \in \mathfrak{R}X$ is \emph{strongly contracting} if there exists
\end{definition}

\subsection{Tits alternative}

\begin{thm}\label{thm:Tits}
Let $G$ be a group acting on a CAT(0) cube complex $X$ of finite dimension~$d$. 
\begin{itemize}
	\item If $G$ does not have a finite orbit in $\overline{X}$, then it contains a non-abelian free subgroup.
	\item Otherwise, $G$ contains a finite-index subgroup which is (locally $X$-elliptic)-by-(free abelian of rank $\leq \dim(X)$).
\end{itemize}
\end{thm}

\begin{proof}
If $G$ has a bounded orbit in $X$, then there is nothing to prove, so, from now on, we assume that $G$ has unbounded orbits in $X$. 

\medskip \noindent
According to \cite[Theorem F]{CapraceSageev}, if $G$ does not have a finite orbit in the visual boundary of $X$, then $G$ contains a non-abelian free subgroup. Otherwise, if $G$ has a finite orbit in the visual boundary of $X$, it follows from \cite[Proposition 2.26]{CFI} that either $G$ has a finite orbit in $\mathfrak{R}X$ or there exists a finite-index subgroup $G' \leq G$ and convex subcomplex $Y \subset \mathfrak{R}X$ such that $G'$ acts on $Y$ with unbounded orbits and with no finite orbit in the visual boundary of $Y$. In the latter case, applying \cite[Theorem F]{CapraceSageev} once again implies that $G$ contains a non-abelian free subgroup; and in the former case, the desired conclusion follows from Theorem \ref{thm:Busemann}.
\end{proof}

\begin{cor}\label{cor:Parabolic}
If a group acts properly on a finite-dimensional CAT(0) cube complex $X$ and fixes a point in $\mathfrak{R}X$, then it must be (locally finite)-by-(free abelian). Conversely, if a (locally finite)-by-(free abelian) group acts on a finite-dimensional CAT(0) cube complex $X$, then it must have a finite orbit in $\mathfrak{R}X$. 
\end{cor}

\section{Proof of the general criterion}

\noindent
This section is dedicated to the proof of the main result of the article, namely:

\begin{thm}\label{thm:MainCriterion}
Let $\mathcal{C}$ be a collection of finitely generated groups. Assume that:
\begin{itemize}
	\item for every group $G \in \mathcal{C}$ and every finite-index subgroup $H \leq G$, there exists a subgroup of $H$ which belongs to $\mathcal{C}$;
	\item for every group $G \in \mathcal{C}$, there exists a subgroup of the commutator subgroup $[G,G]$ which belongs to $\mathcal{C}$;
	\item every group of $\mathcal{C}$ contains a non-abelian free subgroup and a $\delta$-normal subgroup which is (locally finite)-by-(free abelian).
\end{itemize}
Then no group of $\mathcal{C}$ acts properly on a finite-dimensional CAT(0) cube complex. 
\end{thm}

\noindent
Recall that, given a group $G$ and a subgroup $H \leq G$, one says that $H$ is \emph{$\delta$-normal} if, for every finite collection of elements $g_1, \ldots, g_n \in G$, the intersection $\bigcap\limits_{i=1}^n g_iHg_i^{-1}$ is infinite.

\medskip \noindent
For clarity, we decompose the proof into three steps.

\subsection{Step 1: A dichotomy}

\noindent
The first step towards the proof of Theorem \ref{thm:MainCriterion} is to show that, if a group $G$ acts on a finite-dimensional CAT(0) cube complex $X$, then either $G$ stabilises a cubical component of $\mathfrak{R}X$ or $X$ essentially decomposes as a product $X_1 \times \cdots \times X_n$ such that (a finite-index subgroup of) $G$ acts on each factor with a strongly contracting isometry. Such a dichotomy is essentially contained into \cite{CapraceSageev}. 

\begin{prop}\label{prop:dichotomy}
Let $G$ be a group acting on a finite-dimensional CAT(0) cube complex $X$. Assume that $G$ does not stabilise a cubical component of $\mathfrak{R}X$. Then $G$ contains a subgroup $G'$ of index $\leq \dim(X)!$ and $X$ contains a convex subcomplex $Y$ which decomposes as a Cartesian product $X_1 \times \cdots \times X_n$ of $n \geq 1$ irreducible and unbounded cube complexes so that:
\begin{itemize}
	\item $Y$ is $G$-invariant and $G$ preserves the product structure of $Y$;
	\item for every $1 \leq i \leq n$, at least one element of $G'$ defines a strongly contracting isometry of $X_i$.
\end{itemize}
\end{prop}

\noindent
We begin by proving the following preliminary lemma:

\begin{lemma}\label{lem:Fboundary}
Let $G$ be a group acting on a CAT(0) cube complex $X$. If $G$ fixes a point in the visual boundary of $X$, then it stabilises a cubical component of $\mathfrak{R}X$.
\end{lemma}

\begin{proof}
For every CAT(0)-geodesic ray $\sigma$, denote by $\alpha(\sigma)$ the collection of the halfspaces of $X$ in which $\sigma$ eventually lies. The CAT(0)-convexity of halfspaces ensures that $\alpha(\sigma)$ defines an orientation. Moreover, as $\sigma$ goes to infinity, $\alpha(\sigma)$ cannot be principal, so that it defines a point the Roller boundary $\mathfrak{R} X$.

\medskip \noindent
Let $\rho$ be a CAT(0)-geodesic ray representing a point of the visual boundary of $X$ which is fixed by $G$. So, for every $g \in G$, the Hausdorff distance between $\rho$ and $g \cdot \rho$ is finite. Consequently, the number of hyperplanes separating these two rays eventually must be finite, so that the symmetric difference between $\alpha(\rho)$ and $\alpha(g \cdot \rho)=g \cdot \alpha(\rho)$ must be finite. We conclude that $G$ stabilises the cubical component of $\mathfrak{R} X$ which contains $\alpha(\rho)$.
\end{proof}

\begin{proof}[Proof of Proposition \ref{prop:dichotomy}.]
As a consequence of Lemma \ref{lem:Fboundary}, $G$ does not fix a point in the visual boundary of $X$, so that it follows from \cite[Proposition 3.5]{CapraceSageev} that there exists a convex subcomplex $Y \subset X$ on which $G$ acts \emph{essentially}, i.e., no orbit stays within finite Hausdorff distance from a halfspace. Following \cite[Proposition 2.6]{CapraceSageev}, decompose $Y$ as a product of irreducible cube complexes $X_1 \times \cdots \times X_n$, where $n \geq 1$, and let $G'$ be a finite-index subgroup of $G$ which preserves the product structure of $Y$. Notice that, because the action $G \curvearrowright Y$ is essential and fixed-point free in the visual boundary, then so are the actions $G' \curvearrowright X_1, \ldots, X_n$. As a consequence, for every $1 \leq i \leq n$, the cube complex $X_i$ must be unbounded, and \cite[Theorem 6.3]{CapraceSageev} applies so that $G'$ acts on $X_i$ with at least one strongly contracting isometry. 
\end{proof}

\subsection{Step 2: Elementary actions}

\noindent
As a consequence of Proposition \ref{prop:dichotomy}, we need to consider two cases: when the action of our group stabilises a cubical component of the Roller boundary, and when the cube complex decomposes as a product in such a way that (a finite-index subgroup of) the group acts on each factor with a strongly contracting isometry. The section is dedicated to the former case. More precisely, we prove:

\begin{prop}\label{prop:ElementaryAction}
Let $\mathcal{C}$ be a collection of groups satisfying the following two conditions:
\begin{itemize}
	\item for every group $G \in \mathcal{C}$ and every finite-index subgroup $H \leq G$, there exists a subgroup of $H$ which belongs to $\mathcal{C}$;
	\item for every group $G \in \mathcal{C}$, there exists a subgroup of the commutator subgroup $[G,G]$ which belongs to $\mathcal{C}$;
	\item every group of $\mathcal{C}$ is finitely generated and it contains a non-abelian free subgroup.
\end{itemize}
Assume that at least one group of $\mathcal{C}$ acts properly on a finite-dimensional CAT(0) cube complex, and let $d$ denote the smallest dimension of a CAT(0) cube complex on which a group of $\mathcal{C}$ may act properly. If a group $G \in \mathcal{C}$ acts properly on a CAT(0) cube complex $X$ of dimension $d$, then no finite-index subgroup of $G$ stabilises a cubical component of~$\mathfrak{R}X$. 
\end{prop}

\begin{proof}
Assume for contradiction that there exists a group $G \in \mathcal{G}$ which acts properly on a CAT(0) cube complex $X$ of dimension $d$ and which contains a finite-index subgroup stabilising a cubical component $Y \subset \mathfrak{R}X$. In fact, as we know by assumption that such a finite-index subgroup must contain a subgroup which belongs to $\mathcal{C}$, we will suppose without loss of generality that $G$ itself stabilises $Y$. 

\medskip \noindent
Let $\mathcal{H}_1$ denote the collection of the hyperplanes of $X$ which separate at least two vertices of $Y$, and let $\mathcal{H}_2$ denote the complement of $\mathcal{H}_1$. Notice that the image of $Y$ in the Roller boundary of the cubical quotient $X/ \mathcal{H}_1$ is reduced to a single vertex $\xi$. Therefore, $G$ acts on $X/ \mathcal{H}_1$ by fixing the point $\xi$. The Busemann morphism of $\xi$ is only defined on a finite-index of $G$. By assumption, such a subgroup must contain a group $G' \in \mathcal{C}$. Let $\beta : G' \to A$ denote the restriction of the Busemann morphism of $\xi$ to $G'$. Once again, we know that $\mathrm{ker}(\beta)$ has to contain a group $G'' \in \mathcal{C}$. Notice that $G''$ is $X/\mathcal{H}_1$-elliptic as a consequence of Theorem \ref{thm:Busemann}. 

\medskip \noindent
We claim that $G''$ acts properly on $X/ \mathcal{H}_2$. 

\medskip \noindent
For convenience, for $i=1,2$ let $\pi_i$ denote the projection $X \to X/ \mathcal{H}_i$. Fix a vertex $x \in X$. Then, if $S$ denotes the stabiliser of $\pi_2(x)$ in $G''$, we have
$$\begin{array}{lcl} d_X(x,g \cdot x) & = & d_{X/ \mathcal{H}_1}( \pi_1(x), g \cdot \pi_1(x)) + d_{X/ \mathcal{H}_2} (\pi_2(x), g \cdot \pi_2(x)) \\ \\ & = & d_{X/ \mathcal{H}_1} (\pi_1(x), g \cdot \pi_1(x)) \end{array}$$
for every $g \in S$. As $G''$ is $X/ \mathcal{H}_1$-elliptic, it follows that $S$ must be $X$-elliptic. But $G$ acts properly on $X$, which implies that $S$ must be finite, concluding the proof of our claim.

\medskip \noindent
Finally, notice that the dimension of $X/ \mathcal{H}_2$ coincides with the dimension of $Y$, which is smaller than the dimension $d$ of $X$. Thus, we have found a group $G'' \in \mathcal{C}$ which acts properly on a CAT(0) cube complex of dimension $<d$, contradicting the definition of~$d$.
\end{proof}

\subsection{Step 3: Actions with rank-one isometries}

\noindent
We now turn to the second case of the dichotomy provided by Proposition \ref{prop:dichotomy}. Namely, we are interested in actions on cube complexes with strongly contracting isometries. Our main result is the following:

\begin{prop}\label{prop:ContractingNormal}
Let $G$ be a group acting essentially on a CAT(0) cube complex $X$. Assume that no finite-index subgroup of $G$ stabilises a cubical component of $\mathfrak{R} X$ and that $G$ contains at least one strongly contracting isometry. If $G$ contains a subgroup $N$ which is (locally finite)-by-(free abelian), then there exist a finite-index subgroup $N' \leq N$ and elements $a,b,c \in G$ such that the intersection
$$aN'a^{-1} \cap bN'b^{-1} \cap cN'c^{-1}$$
is $X$-elliptic. As a consequence, if $N$ is $\delta$-normal, then $G$ does not act properly on~$X$. 
\end{prop}

\noindent
Recall that a group acts \emph{essentially} on a CAT(0) cube complex if no orbit stays within finite Hausdorff distance from a halfspace.

\medskip \noindent
It is worth noticing that, in the statement of Proposition \ref{prop:ContractingNormal}, the fact that $G$ does not act properly on $X$ when $N$ is $\delta$-normal also follows from that the facts that a group which is not virtually cyclic and which acts on a geodesic metric space with a contracting isometry must be acylindrically hyperbolic \cite{SistoContracting, BBF}, and that s-normal subgroups (which include $\delta$-normal subgroups) of acylindrically hyperbolic groups must be acylindrically hyperbolic as well \cite[Lemma 7.2]{OsinAcyl}. However, Proposition \ref{prop:ContractingNormal} does not only shows that the action is not proper, it identifies an infinite subgroup which has a bounded orbit. Knowing explicitly such a subgroup will be fundamental in the proof of Theorem~\ref{thm:MainCriterion}.

\medskip \noindent
Before turning to the proof of Proposition \ref{prop:ContractingNormal}, we begin by proving two preliminary lemmas.

\begin{lemma}\label{lem:MedianInside}
Let $X$ be a CAT(0) cube complex and $A,B,C$ a facing triple of pairwise strongly separated hyperplanes. Let $A^+,B^+,C^+$ denote respectively the halfspaces delimited by $A,B,C$ which do not contain two hyperplanes among $\{A,B,C\}$. If $\alpha \in A^+$, $\beta \in B^+$ and $\gamma \in C^+$ are three points of $\mathfrak{R}X$, then the median point of $\{\alpha, \beta, \gamma\}$ belongs to $X$.
\end{lemma}

\noindent
Recall that a \emph{facing triple} is the data of three hyperplanes such that no one separates the other two.

\begin{proof}[Proof of Lemma \ref{lem:MedianInside}.]
Fix a vertex $x \in (A^+ \cup B^+ \cup C^+)^c$. Notice that, if $D$ is a halfspace which contains $\alpha$ and $\beta$ but not $x$, then, because $A$ and $B$ are strongly separated, necessarily $D$ has to contain either $A^+$ or $B^+$. Consequently, $D$ separates $x$ from either $A$ or $B$. Similarly, if $D$ contains $\alpha$ and $\gamma$ (resp. $\beta$ and $\gamma)$ but not $x$, then $D$ has to separate $x$ from either $A$ and $C$ (resp. $B$ and $C$). Therefore, if $x$ and the median point $\mu$ of $\{ \alpha, \beta,\gamma\}$ are thought of as orientations, then the symmetric difference between $x$ and $\mu$ has cardinality at most $2(d(x,A)+d(x,B)+d(x,C))$. As a consequence, $\mu$ belongs to the same cubical component as $x$, namely $X$.
\end{proof}

\begin{lemma}\label{lem:DSL}
Let $G$ be a group acting essentially on a CAT(0) cube complex $X$ without fixed point in the visual boundary. Let $A\subset B$ be two halfspaces. For every point $\alpha \in \mathfrak{R}X$ which belongs to $B$, there exists some $g \in G$ such that $g \cdot \alpha \in A$. 
\end{lemma}

\begin{proof}
The Double Skewer Lemma \cite{CapraceSageev} precisely says that there exists some $g \in G$ such that $g \cdot B \subsetneq A$. The desired conclusion follows immediately.
\end{proof}

\noindent
We are now ready to prove our proposition.

\begin{proof}[Proof o Proposition \ref{prop:ContractingNormal}.]
As a consequence of Corollary \ref{cor:Parabolic}, $N$ must have a finite orbit in $\overline{X}$. If $N$ does not have a finite orbit in $\mathfrak{R}X$, then it must be $X$-elliptic, and there is nothing to prove. So suppose that $N$ has a finite orbit in $\mathfrak{R}X$. Up to replacing $N$ with a finite-index subgroup, we suppose without loss of generality that $N$ fixes a point $\xi \in \mathfrak{R}X$. 

\begin{claim}
$X$ contains a facing quadruple of pairwise strongly separated hyperplanes, i.e., four pairwise strongly separated hyperplanes such that no one separates two of the other three.
\end{claim}

\noindent
As $G$ acts on $X$ with a strongly contracting isometry, it follows from Lemma \ref{lem:ContractingPoint} that $\mathfrak{R}X$ contains a strongly contracting point. In fact, because $G$ cannot have a finite orbit in $\mathfrak{R}X$, we know that $\mathfrak{R}X$ has to contain infinitely many strongly contracting points. Let $\alpha, \beta, \gamma, \delta \in \mathfrak{R} X$ be four of them. Fix a decreasing sequence of halfspaces $(A_i)$ (resp. $(B_i)$, $(C_i)$, $(D_i)$) bounded by pairwise strongly separated hyperplanes and such that $\alpha$ (resp. $\beta$, $\gamma$, $\delta$) belongs to $A_i$ (resp. $B_i$, $C_i$, $D_i$) for every $i$. 

\medskip \noindent
Notice that, if $A_i$ and $B_j$ intersect for some $i,j$, then there must exist $p,q$ such that $A_p \subset B_j$ and $B_q \subset A_i$. As $\alpha$ and $\beta$ are distinct, we deduce from Lemma \ref{lem:ContractingPoint} that there exists some $r$ such that $A_i$ and $B_j$ are disjoint for every $i,j \geq r$. A similar conclusion holds for any pair of sequences among $(A_i)$, $(B_i)$, $(C_i)$, $(D_i)$. As a consequence, there exists some $s$ such that $A_i$, $B_j$, $C_k$ and $D_\ell$ are pairwise disjoint for every $i,j,k,\ell \geq s$. 

\medskip \noindent
Therefore, the hyperplanes bounding the four halfspaces $A_r$, $B_r$, $C_r$ and $D_r$ define a facing quadruple. It follows that the hyperplanes bounding $A_{r+1}$, $B_{r+1}$, $C_{r+1}$ and $D_{r+1}$ define a facing quadruple of pairwise strongly separated hyperplanes, concluding the proof of our claim. 

\medskip \noindent
Let $A,B,C,D$ be a facing quadruple of pairwise strongly separated hyperplanes, and let $A^+,B^+,C^+,D^+$ denote the halfspace delimited respectively by $A,B,C,D$ which are disjoint from the three other hyperplanes. Of course, the point $\xi$ of $\mathfrak{R}X$ which is fixed by $N$ may belong to at most one of these four halfspace. Up to relabelling our hyperplane, let us say that $\xi$ does not belong $A^+$, $B^+$ and $C^+$. By applying Lemma \ref{lem:DSL} three times, we find three elements $a,b,c \in G$ such that $a \xi \in A^+$, $b \xi \in B^+$ and $c \xi \in C^+$. The median point $\mu$ of these three points of $\overline{X}$, which is a vertex of $X$ according to Lemma \ref{lem:MedianInside}, is stabilised by
$$aNa^{-1} \cap bNb^{-1} \cap cNc^{-1},$$
concluding the proof of the theorem.
\end{proof}

\subsection{Conclusion}

\noindent
By combining Propositions \ref{prop:dichotomy}, \ref{prop:ElementaryAction} and \ref{prop:ContractingNormal}, we are now ready to prove the main criterion provided by Theorem \ref{thm:MainCriterion}.

\begin{proof}[Proof of Theorem \ref{thm:MainCriterion}.]
Assume for contradiction that there exists a group of $\mathcal{C}$ which acts properly on a finite-dimensional CAT(0) cube complex, and let $d$ denote the smallest dimension of a CAT(0) cube complex on which a group of $\mathcal{C}$ may act properly. Now fix a group $G \in \mathcal{C}$ acting properly on a CAT(0) cube complex of dimension $d$. According to Proposition \ref{prop:ElementaryAction}, no finite-index subgroup $G$ stabilises a cubical component of $\mathfrak{R}X$. It follows from Proposition \ref{prop:dichotomy}, that $G$ contains a finite-index subgroup $G'$ and $X$ contains a convex subcomplex $Y$ which decomposes as a Cartesian product $X_1 \times \cdots \times X_n$ of $n \geq 1$ irreducible and unbounded cube complexes so that:
\begin{itemize}
	\item $Y$ is $G$-invariant and $G$ preserves the product structure of $Y$;
	\item for every $1 \leq i \leq n$, at least one element of $G'$ defines a strongly contracting isometry of $X_i$.
\end{itemize}
Notice that, if $N$ is an $\delta$-normal subgroup of $G$ which is (locally finite)-by-(free abelian), then $N':=N \cap G'$ defines a similar subgroup for $G'$. 

\medskip \noindent
Next, we deduce from Proposition \ref{prop:ContractingNormal} that $N'$ contains a finite-index subgroup $N''$ such that, for every $1 \leq i \leq n$, there exist elements $a_i,b_i,c_i \in G'$ so that the intersection
$$a_iN''a_i^{-1} \cap b_iN''b_i^{-1} \cap c_iN''c_i^{-1}$$
is $X_i$-elliptic. Consequently, the intersection
$$E:= \bigcap\limits_{i=1}^n \left( a_iN''a_i^{-1} \cap b_iN''b_i^{-1} \cap c_iN''c_i^{-1} \right)$$
must be $X_i$-elliptic for every $1 \leq i \leq n$. In other words, $E$ must be $Y$-elliptic. Since $N$ (and a fortiori $N''$ as a finite-index subgroup of $N$) is $\delta$-normal, we know that $E$ is infinite, and we conclude that $G$ does not act properly on $X$. 
\end{proof}

\section{Applications}

\subsection{Wreath products}

\noindent
By applying Theorem \ref{thm:MainCriterion} to the specific case where our collection $\mathcal{C}$ is reduced to a single finitely generated group, one obtains immediately the following statement:

\begin{cor}
Let $G$ be a finitely generated group. Assume that:
\begin{itemize}
	\item the commutator subgroup $[G,G]$ and every finite-index subgroup of $G$ contains a copy of $G$;
	\item $G$ contains a non-abelian free subgroup and a $\delta$-normal subgroup which is (locally finite)-by-(free abelian).
\end{itemize}
Then $G$ does not act properly on a finite-dimensional CAT(0) cube complex.
\end{cor}

\noindent
A straightforward application of this corollary is the following proposition:

\begin{prop}\label{prop:Wreath}
For every non-trivial finite group $F$, the wreath product $F \wr \mathbb{F}_2$ does not act properly on a finite-dimensional CAT(0) cube complex.
\end{prop}

\noindent
By combining this proposition with Tits alternative for finite-dimensional CAT(0) cube complexes, one obtains:

\begin{cor}\label{cor:WreathFiniteDim}
Let $A$ and $B$ be two non-trivial groups. If the wreath product $A \wr B$ acts properly on a finite-dimensional CAT(0) cube complex, then $A$ must be a torsion group and $B$ must be virtually (locally finite)-by-(free abelian of finite rank). 
\end{cor}

\begin{proof}
If $A$ contains an infinite-order element, then $A \wr B$ has to contain a free abelian group of infinite rank. Therefore, $A \wr B$ cannot act properly on a finite-dimensional CAT(0) cube complex in this case. Next, it follows from Proposition \ref{prop:Wreath} that, if $A \wr B$ acts properly on a finite-dimensional CAT(0) cube complex, then $B$ cannot contain a non-abelian free subgroup. We conclude from Theorem \ref{thm:Tits} that $B$ must be virtually (locally finite)-by-(free abelian). 
\end{proof}

\noindent
We do not know whether or not the converse of Corollary \ref{cor:WreathFiniteDim} holds. Consequently, the classification of wreath products acting properly on finite-dimensional CAT(0) cube complexes remains unknown.

\subsection{Normalisers of infinite torsion subgroups}

\noindent
In this section, our goal is to deduce from Theorem \ref{thm:MainCriterion} some knowledge about normalisers of infinite torsion subgroups in groups acting properly on finite-dimensional CAT(0) cube complexes. Roughly speaking, we show that, if a group $G$ contains an infinite torsion subgroup $L$ and if $G$ acts properly on a finite-dimensional CAT(0) cube complex, then either the normaliser $N_G(L)$ is close to be free abelian or, for every $k \geq 1$, $N_G(L)$ contains a non-abelian free subgroup commuting with a subgroup of $L$ of size $\geq k$. More precisely, the main result of this section is:

\begin{thm}\label{thm:TorsionNormaliser}
Let $G$ be a group acting properly on a finite-dimensional CAT(0) cube complex $X$. Assume that $G$ contains an infinite torsion subgroup $L$. Either the normaliser $N_G(L)$ contains a finite-index subgroup which is (locally finite)-by-$\mathbb{Z}^n$ for some $n \leq \dim(X)$; or, for every $k \geq 1$, $N_G(L)$ a contains a non-abelian free subgroup centralising a subgroup of $L$ of cardinality $\geq k$.
\end{thm}

\noindent
By begin by proving the following preliminary statement:

\begin{prop}\label{prop:TorsionNormal}
Let $G$ be a finitely generated group which contains a non-abelian free subgroup $F$ and a normal subgroup $L$ which is infinite torsion. If $G$ acts properly on a finite-dimensional CAT(0) cube complex, then, for every $k \geq 1$, $F$ must contain a non-abelian subgroup which centralises a subgroup of $L$ of cardinality $\geq k$.
\end{prop}

\begin{proof}
Let $\mathcal{C}$ denote the collection of finitely generated groups $G$ such that:
\begin{itemize}
	\item $G$ contains a non-abelian free subgroup $F$ and a normal subgroup $L$ which is infinite torsion;
	\item there exists some $k \geq 1$ such that no subgroup of $L$ of cardinality $\geq k$ is centralised by a non-abelian subgroup of $F$.
\end{itemize}
In order to prove our proposition, it is sufficient to show that $\mathcal{C}$ satisfies the conditions of Theorem \ref{thm:MainCriterion}. 

\medskip \noindent
So fix a group $G \in \mathcal{G}$. If $H$ is a finite-index subgroup of $G$, it is clear that $F \cap H$ and $L \cap H$ are subgroups of $H$ allowing us to deduce that $H$ has to belong to $\mathcal{C}$. 

\medskip \noindent
Now, assume that $H$ is the commutator subgroup of $G$. As $G$ is finitely generated, the quotient $G/H$ must be a finitely generated abelian group. Since the torsion of such an abelian group is bounded, it follows that $L \cap H$ must have finite index in $L$; a fortiori, it is infinite. We also know that $F \cap H$ must be non-abelian as it contains $[F,F]$. However, $H$ may not be finitely generated so we cannot deduce that $H$ belongs to $\mathcal{C}$. Fix a subset $S \subset L \cap H$ of cardinality more than $k$ (the constant associated to $G$ by the definition of $\mathcal{C}$) and a finitely generated non-abelian subgroup $F' \leq F \cap H$. We claim that $\langle S, F' \rangle$ belongs to $\mathcal{C}$. Of course, $\langle S, F' \rangle$ is finitely generated by construction, and it contains $F'$ as a non-abelian free subgroup. Let $L'$ denote the normal closure of $S$ in $\langle S, F' \rangle$. If $L'$ is finite, the kernel of the action of $F'$ on $L'$ by conjugations defines a non-abelian subgroup of $F$ centralising the subgroup $\langle S \rangle$ of $L$, which has cardinality more than $k$. By definition of $k$, this is impossible, so $L'$ must be infinite. Moreover, by definition of $k$, no subgroup of $L' \leq L$ of cardinality $ \geq k$ can be centralised by a non-abelian subgroup of $F' \leq F$. We conclude that indeed $\langle S,F' \rangle \leq H$ belongs to $\mathcal{C}$.
\end{proof}

\begin{proof}[Proof of Theorem \ref{thm:TorsionNormaliser}.]
If the normaliser $N_G(L)$ does not contain a non-abelian free subgroup, it follows from Theorem \ref{thm:Tits} that it has be virtually (locally finite)-by-(free abelian of rank $\leq \dim(X)$), and there is nothing else to prove. From now on, assume that $N_G(L)$ contains a non-abelian finitely generated free subgroup $F$. For every $k \geq 1$, fix a finite subgroup $S_k \leq L$ of cardinality $\geq k$. Notice that $\langle S_k, F \rangle$ is finitely generated for every $k$. We distinguish two cases.

\medskip \noindent
If, for every $k \geq 1$, the normal closure $L_k$ of $S_k$ in $\langle L_k,F \rangle$ is finite, then the kernel of the action of $F$ on $L_k$ by conjugations defines a non-abelian free subgroup of $F \leq N_G(L)$ centralising $L_k$, which is a subgroup of $L$ of cardinality $\geq k$. The desired conclusion follows.

\medskip \noindent
Otherwise, if there exists some $k \geq 1$ such that $L_k$ is infinite, then $\langle S_k, F \rangle$ defines a finitely generated group which contains $F$ as a non-abelian free subgroup, which contains $L_k$ as a normal infinite torsion subgroup, and which acts properly on a finite-dimensional CAT(0) cube complex. The desired conclusion follows from Proposition \ref{prop:TorsionNormal}.
\end{proof}

\addcontentsline{toc}{section}{References}

\bibliographystyle{alpha}
{\footnotesize\bibliography{TorsionCC}}

\end{document}